\documentclass[a4paper, 12pt, reqno]{amsart}
\usepackage[left=2.5cm,right=2.5cm,top=3cm,bottom=3cm,a4paper]{geometry}
\usepackage{amsfonts}
\usepackage{amssymb}
\usepackage{amsmath}
\usepackage[doublespacing]{setspace}

\newtheorem{lem}{Lemma}
\newtheorem{thm}{Theorem}

\title{New simultaneous methods for finding all zeros of a polynomial}

\begin{document}
\author[Jun-Seop Song]{Jun-Seop Song$^{a,*}$\\ $^a$College of Medicine, Yonsei University, Seoul, Republic of Korea}
\thanks{*Corresponding author}
\email{junseopsong@yonsei.ac.kr}



\maketitle

\begin{abstract}
The purpose of this paper is to present three new methods for
finding all simple zeros of polynomials simultaneously. First, we
give a new method for finding simultaneously all simple zeros of
polynomials constructed by applying the Weierstrass method to the
zero in the trapezoidal Newton's method, and prove the convergence
of the method. We also present two modified Newton's methods
combined with the derivative-free method, which are constructed by
applying the derivative-free method to the zero in the trapezoidal
Newton's method and the midpoint Newton's method, respectively.
Finally, we give a numerical comparison between various simultaneous
methods for finding zeros of a polynomial.
\end{abstract}

\section{Introduction}
With a typical iteration method such as Newton's method, an initial
approximation of a zero converges to a specific zero, but the
Weierstrass method (or Durand-Kerner method) approximates all simple
(real or complex) zeros of polynomial simultaneously (see
\cite{Durand, Kerner}).

Let $P(z)=z^n+a_1z^{n-1}+\cdots+a_{n-1}z+a_n$ be a polynomial of
degree $n$ having simple zeros with constants $a_1, a_2, \ldots,
a_n$. Let $r_1, \ldots, r_n$ be the distinct zeros of $P(z)$ and let
distinct complex numbers $z_1, \ldots, z_n$ be their approximations.
The {\it Weierstrass method} (Durand-Kerner method) is defined as
\begin{equation} \label{WDK-fullform}
z_i^{m+1} = z_i^{m} - \frac{P(z_i^m)}{\prod_{j \neq
i}(z_i^m-z_j^m)},
\end{equation}
for $m \ge 0$, and this method is one of the most frequently used
iterative methods which give simultaneous computation of all zeros
of $P$. If a function $W_i(z)$ is defined by
$$W_i(z) = \frac{P(z)}{\prod_{j \neq i}(z-z_j)},$$
then $W_i(z)$ has the same zeros as the polynomial $P$, and so the
problem of finding the zeros of $P$ reduces to that of zeros of the
function $W_i (z)$. If we denote $W_i = W_i(z_i)$ for $i=1, 2,
\ldots, n$ in the case of $z=z_i$, (\ref{WDK-fullform}) can be
written as
\begin{equation} \label{WDK}
\hat{z_i} = z_i - W_i,
\end{equation}
where $z_i$ is a current approximation and $\hat{z_i}$ is a new
approximation to a zero of polynomial $P(z)$. The method constructed
by (\ref{WDK}) is called the {\it Weierstrass-like method} (briefly, WLM).

The aim of this paper is to present three new methods for finding
all simple zeros of polynomials simultaneously. These new methods
are based on the Frontini-Sormani's midpoint Newton's method
(\cite{Traub}) and the Weerakoon's trapezoidal Newton's method
(\cite{VNM}) which were modifications of the Newton's method through
iterative approximations.

It is well known that Newton's method is defined by
$x^*=x-\frac{f(x)}{f'(x)}$ with an approximation $x$ and a new
approximation $x^*$ of a zero, and is efficient to find a zero of an
equation $f(x)=0$ for a differentiable function $f$ with proper
conditions and a sufficiently close initial value (see \cite{VNM}).

In \cite{VNM}, Weerakoon proposed the {\it trapezoidal Newton's
method} defined by
\begin{equation} \label{trapezoidal-method}
\hat{x}=x-\frac{2f(x)}{f'(x)+f'(x^*)}.
\end{equation}
He applied Newton's method to the $x^*$ of the denominator.

Along with (\ref{trapezoidal-method}), the {\it midpoint Newton's
method} that Frontini-Sormani proposed in \cite{Traub} is
constructed as
\begin{equation} \label{midpoint-method}
\hat{x} = x - \frac{f(x)}{f'\left(x+\frac{1}{2}(x^*-x)\right)}.
\end{equation}
They also applied Newton's method to the $x^*$ of the denominator,
and so set $x^*=x-\frac{f(x)}{f'(x)}$.

Both the trapezoidal Newton's method and the midpoint Newton's
method are of cubic order, while the original Newton's method was of
quadratic order. A variety of methods can be applied to the $x^*$ in
addition to Newton's method. Petkovi\'c \cite{Petkovic-midpoint}
derived the following simultaneous method for finding all simple
zeros of polynomials by applying the Weierstrass method to the $x^*$
in the midpoint Newton's method:
\begin{equation} \label{pm}
\hat{z_i} = z_i - \frac{P(z_i)}{P'\left(z_i-\frac{1}{2}W_i\right)},
\end{equation}
which is called {\it Newton-Weierstrss method} (or NWM). Also,
Petkovi\'c \cite{Petkovic-derivative-free} found the following {\it
derivative-free method} (or DFM) defined as
\begin{equation} \label{derivative-free}
\hat{z_i} = z_i - \frac{W_i}{1-P(z_i-W_i)/P(z_i)},
\end{equation}
which has a similar form with the one above and this method is of
cubic order.

In this paper, we present three new methods for the simultaneous
approximation of all simple zeros of polynomials by applying the
Weierstrass-like method and the derivative-free method to $x^*$ in
the trapezoidal Newton's method and the midpoint Newton's method.

Throughout this paper, the convergence of zeros will be discussed
and the order will be calculated for new constructed methods. We
will use the notation $a=O_M(b)$ for two complex numbers $a$ and $b$
whose moduli are of the same order, that is,
$|a|=O\left(|b|\right)$. In addition, the error is defined as $|e| =
\underset{i=1,\ldots,n}{\max} \{ |e_i| \}$ with $e_i=z_i-r_i$ for
$i=1,\ldots,n$.

In all discussions, the order related to $e_i$, which is an error of
the previously approximated zeros $z_i$, is presumed to be the same.
After that, we will show that the order related to $e_i$, which is
an error of the approximated zeros concerning each method, is
identical. For the same being, the order related to the already
approximated zeros $\hat{e_i}$ is hypothesized to be identical as
follows:
$$e_i = O_M(e)~~ \text{for all } i.$$

In Section 2, we give a new method for finding simultaneously all
simple zeros of polynomials constructed by applying the Weierstrass
method to the $x^*$ in the trapezoidal Newton's method, and prove
the convergence of the method. In Section 3, we present two modified
Newton's methods combined with the derivative-free method. They are
constructed by applying the derivative-free method to the $x^*$ in
the trapezoidal Newton's method and the midpoint Newton's method,
respectively. In Section 4, we give a numerical comparison between
various simultaneous methods for finding zeros of a polynomial.
Finally, we conclude that the convergence of all new constructed
methods in this paper are similar or superior than other iterative
methods of cubic order.

\section{Weierstrass-like Trapezoidal Newton's method}

In this section, we construct a new method for finding
simultaneously all simple zeros of polynomials of cubic order. By
applying the Weierstrass method (\ref{pm}) to the $x^*$ in the
trapezoidal Newton's method (\ref{trapezoidal-method}), we derive a
new method constructed as follows:
\begin{equation} \label{method1}
\hat{z_i}=z_i-\frac{2P(z_i)}{P'(z_i)+P'\left(z_i-W_i\right)}.
\end{equation}
We call (\ref{method1}) the {\it Weierstrass-like trapezoidal Newton's method},
and from this, simply, call it {\bf Method 1}.

The calculation and discussion of the order of Method 1 are similar
to those of the Newton-Weierstrass method, which is an alteration of
Petkovi\'c's midpoint Newton's method (see
\cite{Petkovic-midpoint}). From (\ref{method1}), we have the
following theorem.

\begin{lem} \label{lemma1}
For a polynomial $P(z)$, we have $$
\frac{P(z_i)}{P'(z_i)}\left( \frac{P''(z_i)P(z_i)}{2P'(z_i)^2}O_M(e)
+ O_M(e^2) \right) = O_M(e^3).
$$
\end{lem}

\begin{proof}
By the Taylor's expansion around $r_i$, we have that
\begin{eqnarray}
\label{eq8} &&P(z_i) = P'(r_i)\left(e_i + 
\frac{P''(r_i)}{2P'(r_i)}e_i^2 + O_M(e^3)\right),\notag \\
&&P'(z_i) = P'(r_i)\left( 1 + \frac{P''(r_i)}{P'(r_i)}e_i +
\frac{P'''(r_i)}{2P'(r_i)}e_i^2 + O_M(e^3) \right),\\
&&P''(z_i) = P'(r_i)\left( \frac{P''(r_i)}{P'(r_i)} +
\frac{P'''(r_i)}{P'(r_i)}e_i + \frac{P''''(r_i)}{2P'(r_i)}e_i^2 +
O_M(e^3) \right).\notag
\end{eqnarray}
From (\ref{eq8}), we obtain
$$\aligned &\frac{P(z_i)}{P'(z_i)}\left(
\frac{P''(z_i)P(z_i)}{2P'(z_i)^2}O_M(e)
+ O_M(e^2) \right)  \\
 &= \frac{1}{2}\left( e_i + \frac{P''(r_i)}{2P'(r_i)}e_i^2 + O_M(e^3) \right)^2 \left( 1 - \frac{P''(r_i)}{P'(r_i)}e_i
 + \frac{P'''(r_i)}{2P'(r_i)}e_i^2 + O_M(e^3) \right)^3 \\
& \quad \times \left( \frac{P''(r_i)}{P'(r_i)} +
\frac{P'''(r_i)}{P'(r_i)}e_i + \frac{P''''(r_i)}{2P'(r_i)}e_i^2 +
O_M(e^3) \right)O_M(e) +
O_M(e^3) \\
& = \frac{1}{2}\left(O_M(e)\right)^2 \left(1+O_M(e)\right)^3
\left(\frac{P''(r_i)}{P'(r_i)}+O_M(e)\right)O_M(e) + O_M(e^3) =
O_M(e^3).  \endaligned$$ \qedhere
\end{proof}

From Lemma \ref{lemma1}, we have the following theorem.
\begin{thm}
If the approximate zero $x_i$ grounded from Method 1 is close enough
to $r_i$   and the order of $e_i$ is the same, then the order of
$\hat{e_i}$ is identical, and $|\hat{e_i}| = O_M\left(|e|^3\right)$ is
formed.
\end{thm}

\begin{proof} We easily see that the following equation is satisfied.
$$P(z_i) = \prod_{j=1}^{n} (z_i-r_j) = (z_i-r_i) \prod_{j \neq i}(z_i-r_j) = e_i \prod_{j \neq i}(z_i-r_j) = O_M(e). $$
That is,
\begin{equation} \label{thm1-a}
W_i=W_i(z_i)=O_M\left(P(z_i)\right)=O_M(e).
\end{equation}

If $Q(z)=P(z)-\prod_{j=1}^n(z-z_j)$, then $Q(z)$ is a polynomial of
order  $n-1$, and $Q(z_i)=P(z_i)$ for all $i$. Therefore, $Q(z)$ is
the Lagrange interpolation of points $z_1, z_2, \ldots, z_n$, and so
we have
$$\aligned Q(z) &= \sum_{j=1}^n\left( P(z_j)\prod_{k \neq
j}\frac{z-z_k}{z_j-z_k} \right)=\sum_{j=1}^n\left( \left(
W_j\prod_{k \neq j}(z_j-z_k)
\right)\prod_{k \neq j}\frac{z-z_k}{z_j-z_k} \right)\\
&=\sum_{j=1}^n\left( W_j\prod_{k \neq j} (z-z_k) \right) =
\left(\sum_{j=1}^n
\frac{W_j}{z-z_j}\right)\prod_{k=1}^n(z-z_k).\endaligned$$
Therefore, we obtain
\begin{eqnarray} \label{aaaaa}
P(z) & =& \left(1 + \sum_{j=1}^n
\frac{W_j}{z-z_j}\right)\prod_{j=1}^n(z-z_j) \\
&=& W_i \prod_{j \neq i}(z-z_j) + \left(1 + \sum_{j \neq i}
\frac{W_j}{z-z_j}\right)\prod_{j=1}^n(z-z_j). \notag
\end{eqnarray}
From (\ref{aaaaa}), it follows that
\begin{equation} \label{aaabb}
\frac{P'(z)}{P(z)} = \left(\sum_{j \neq i} \frac{1}{z-z_j}\right) +
\frac{1 + \sum_{j \neq i}\frac{W_j}{z-z_j} - (z-z_i)\sum_{j \neq
i}\frac{W_j}{(z-z_j)^2}}{W_i + (z-z_i)\left(1+\sum_{j \neq
i}\frac{W_j}{z-z_j}\right)}.
\end{equation}
Substituting $z=z_i$ in (\ref{aaabb}), we obtain
$$\frac{P'(z_i)}{P(z_i)} = \left(\sum_{j \neq i} \frac{1}{z_i-z_j}\right) + \frac{1 + \sum_{j \neq i}\frac{W_j}{z_i-z_j}}{W_i}.$$
Therefore, we have
\begin{equation} \label{thm1-b}
W_i = \frac{P(z_i)}{P'(z_i)}\cdot\frac{1 + \sum_{j \neq
i}\frac{W_j}{z_i-z_j}}{1-\frac{P(z_i)}{P'(z_i)}\sum_{j \neq
i}\frac{1}{z_i-z_j}} = \frac{P(z_i)}{P'(z_i)}\left(1 +
O_M(e)\right).
\end{equation}

Now we will find the order of Method 1. If the Taylor's expansion is
applied to $P'(z_i-W_i)$, then we have
$$\aligned \hat{z_i} &=z_i-\frac{2P(z_i)}{P'(z_i)+P'\left(z_i-W_i\right)}\\
&=z_i-\frac{2P(z_i)}{P'(z_i)+\left( P'(z_i) -
P''(z_i)W_i + \frac{1}{2}P'''(z_i)W_i^2 +\cdots\right)}\\
&=z_i-\frac{P(z_i)}{P'(z_i)-\frac{1}{2}P''(z_i)W_i+\frac{1}{4}P'''(z_i)W_i^2+\cdots}
\, .
\endaligned$$
By (\ref{thm1-a}), (\ref{thm1-b}) and Lemma \ref{lemma1}, we have that
$$\aligned \hat{z_i} &=z_i-\frac{P(z_i)}{P'(z_i)\left( 1-\frac{P''(z_i)}{2P'(z_i)}W_i+O_M(e^2)
\right)}\\&=z_i-\frac{P(z_i)}{P'(z_i)}\cdot\left(1+\frac{P''(z_i)}{2P'(z_i)}W_i+O_M(e^2)
\right)\\&= z_i - \frac{P(z_i)}{P'(z_i)}\cdot \left( 1+
\frac{P''(z_i)P(z_i)}{2P'(z_i)^2}\left(1+O_M(e)\right) + O_M(e^2)
\right) \\&= z_i - \frac{P(z_i)}{P'(z_i)} -
\frac{P''(z_i)P(z_i)^2}{2P'(z_i)^3} - \frac{P(z_i)}{P'(z_i)} \left(
\frac{P''(z_i) P(z_i)}{2 P'(z_i)^2} O_M(e) + O_M(e^2) \right)
\\&=z_i - \frac{P(z_i)}{P'(z_i)} -
\frac{P''(z_i)P(z_i)^2}{2P'(z_i)^3}-O_M(e^3).
\endaligned $$
The Chebyshev's method  is defined by
$$\hat{x}=x-\left( 1 + \frac{f''(x)f(x)}{2f'(x)^2} \right)\frac{f(x)}{f'(x)},$$
and of cubic order (see \cite[Section 5.2]{Traub}). According to the
Chebyshev's method, we see that
$$ z_i - \frac{P(z_i)}{P'(z_i)} - \frac{P''(z_i)P(z_i)^2}{2P'(z_i)^3} - r_i = O_M(e^3). $$
Therefore, the order of $\hat{e_i}$ is calculated as follows:
$$ |\hat{e_i}| = |\hat{z_i}-r_i| = \left|z_i - \frac{P(z_i)}{P'(z_i)} - \frac{P''(z_i)P(z_i)^2}{2P'(z_i)^3} - r_i + O_M(e^3)\right|
= O_M\left(|e|^3\right).$$ \qedhere
\end{proof}

\section{Modified Newton's methods combined with Derivative-free method}

In this section, we present two modified Newton's methods combined
with the derivative-free method (\ref{derivative-free}) for finding
all simple zeros of a polynomials simultaneously. The one is a form
that the derivative-free method is applied to the $x^*$ in the
trapezoidal Newton's method (\ref{trapezoidal-method}) as follows:
\begin{equation} \label{method2}
\hat{z_i}=z_i-\frac{2P(z_i)}{P'(z_i) +
P'\left(z_i-\frac{W_i}{1-P(z_i-W_i)/P(z_i)}\right)},
\end{equation}
which is called the {\it Derivative-free trapezoidal Newton's
method}, or simply, {\bf Method 2}.

From (\ref{method2}), we have the following theorem.

\begin{thm}
If the approximate zero $x_i$ grounded from Method 2 is close enough
to $r_i$  and the order of $e_i$ is the same, then the order of
$\hat{e_i}$ is identical, and $|\hat{e_i}| = O_M\left(|e|^3\right)$
is formed.
\end{thm}
\begin{proof}
Since Petkovi\'c's derivative-free method (\ref{derivative-free}) is
of cubic order
\begin{equation} \label{derivative-free-order}
z_i-\frac{W_i}{1-P(z_i-W_i)/P(z_i)}-r_i=O_M(e^3),
\end{equation}
(see \cite{Petkovic-derivative-free}). Using
(\ref{derivative-free-order}) and the Taylor's expansion,
$\hat{z_i}$ is calculated as follows. (In this case, $C_j =
\frac{1}{j!}\cdot\frac{P^{(j)}(r_i)}{P'(r_i)}$).
$$ \aligned
\hat{z_i}&=z_i-\frac{2P(z_i)}{P'(z_i) +
P'\left(z_i-\frac{W_i}{1-P(z_i-W_i)/P(z_i)}\right)}\\
&= z_i - \frac{2
P(r_i+e_i)}{P'(r_i+e_i)+P'\left(r_i+O_M(e^3)\right)}\\ &=z_i -
\frac{2\left(P(r_i)+P'(r_i)e_i+\frac{1}{2}P''(r_i)e_i^2+\frac{1}{6}
P'''(r_i)e_i^3+\cdots\right)}{\left(P'(r_i)+P''(r_i)e_i+\frac{1}{2}P'''(r_i)e_i^2+\cdots\right)+\left(P'(r_i)+P''(r_i)O_M(e^3)+\cdots\right)}\\
&=z_i - \frac{2P'(r_i)\left( e_i+C_2e_i^2+C_3e_i^3+O_M(e^4)
\right)}{2P'(r_i)\left( 1+C_2e_i+\frac{3}{2}C_3e_i^2+O_M(e^3)
\right)}\\ &=z_i-e_i\left(1+C_2e_i+C_3e_i^2+O_M(e^3)\right)\left(
1-C_2e_i+\frac{3}{2}C_3e_i^2+O_M(e^3) \right)\\
&=z_i-e_i+O_M(e^3).\endaligned$$

Therefore, the order of $\hat{e_i}$ is calculated as follows.
$$ |\hat{e_i}| = |\hat{z_i}-r_i| = \left|z_i-e_i+O_M(e^3) - r_i\right| = O_M\left(|e|^3\right).$$
\qedhere
\end{proof}

Now we apply the derivative-free method to the $x^*$ in the midpoint
Newton's method (\ref{midpoint-method}) and construct the iteration
as follows:
\begin{equation} \label{method3}
\hat{z_i} = z_i -
\frac{P(z_i)}{P'\left(z_i-\frac{1}{2}\cdot\frac{W_i}{1-P(z_i-W_i)/P(z_i)}\right)}\,
,
\end{equation}
which is  called the {\it Derivative-free midpoint Newton's method}.
From this, we call it {\bf Method 3} simply. From (\ref{method3}),
we have the following theorem.

\begin{thm}
If the approximate zero $x_i$ grounded from Method 3 is close enough
to $r_i$  and the order of $e_i$ is the same, then the order of
$\hat{e_i}$ is identical, and $|\hat{e_i}| = O_M\left(|e|^3\right)$
is formed.
\end{thm}
\begin{proof}
By using (\ref{derivative-free-order}) and Taylor's expansion,
$\hat{z_i}$ is calculated as follows: (In this case, $C_j =
\frac{1}{j!}\cdot\frac{P^{(j)}(r_i)}{P'(r_i)}$).
$$ \aligned \hat{z_i} &= z_i -
\frac{P(z_i)}{P'\left(z_i-\frac{1}{2}\cdot\frac{W_i}{1-P(z_i-W_i)/P(z_i)}\right)}\\
&= z_i -
\frac{P(z_i)}{P'\left(\frac{1}{2}z_i+\frac{1}{2}r_i+O_M(e^3)\right)}\\
&=z_i-\frac{P(r_i+e_i)}{P'\left( r_i+\frac{1}{2}e_i+O_M(e^3)
\right)}\\
&=z_i-\frac{P(r_i)+P'(r_i)e_i+\frac{1}{2}P''(r_i)e_i^2+\frac{1}{6}P'''(r_i)e_i^3
+\cdots}{P'(r_i)+P''(r_i)\left(\frac{1}{2}e_i+O_M(e^3)\right)+\frac{1}{2}P'''(r_i)\left(\frac{1}{2}e_i+O_M(e^3)\right)^2+\cdots}\\
&=z_i-\frac{P'(r_i)\left( e_i+C_2e_i^2+C_3e_i^3+O_M(e^4)
\right)}{P'(r_i)\left( 1+C_2e_i+\frac{3}{4}C_3e_i^2+O_M(e^3)
\right)}\\ &=z_i-e_i\left( 1+C_2e_i+C_3e_i^2+O_M(e^3)\right)\left(
1-C_2e_i+\frac{3}{4}C_3e_i^2+O_M(e^3) \right)\\
&=z_i-e_i+O_M(e^3).\endaligned $$ Therefore, the order of
$\hat{e_i}$ can be calculated as follows.
$$ |\hat{e_i}| = |\hat{z_i}-r_i| = \left|z_i-e_i+O_M(e^3) - r_i\right| = O_M\left(|e|^3\right).$$
\qedhere
\end{proof}

\section{Numerical comparison}
In this section, we give numerical experiments and comparisons
between various simultaneous methods for finding zeros of a
polynomial. These methods are all of cubic order. They include
Method 1, Method 2, Method 3, the derivative-free method (DFM), the
Petkovi\'c's Newton-Weierstrass method (NWM), and the
Weierstrass-like method (WLM).

For a polynomial $P(z)=z^n+a_1z^{n-1}+\cdots+a_{n-1}z+a_n$, we
choose initial approximations as Aberth's approach (see
\cite{Aberth}):
$$z_k^{(0)} = -\frac{a_1}{n}+R \exp\left(\frac{i\pi}{n}\left(2k-\frac{3}{2}\right)\right).$$
In this case, R is a radius of a circle, where the initial zeros by
Aberth's approach are located in complex number plane. We use the
following Henrici's fomula to select R (see \cite{Henrici}):
$$R=2 \max_{1 \leq k \leq n} |a_k|^{1/k}.$$
According to Henrici's formula, a disk $\{z : |z| < R\}$ centered at
the origin contains all zeros of polynomial $P(z)$.

The polynomials that we used on numerical comparison are as follows:
\begin{eqnarray} \label{test-polynomials}
&&P_1 = (x-1)(x-2)(x-3)(x-4) \nonumber \\
&&P_2 = (x-1)(x-2)(x-3)(x-4)(x-5) \\
&&P_3 = (x-1)(x-2)(x-3)(x-4)(x-5)(x-6) \nonumber \\
&&P_4 = x^8 + 5x^7 + 3x^6 + 7x^5 + 6x^4 + 8x^3 + x^2 + 3x + 7 \nonumber
\end{eqnarray}
Here $P_1 (x), P_2 (x)$ and  $P_3 (x)$ are Wilkinson's polynomials
when $n=4,5,6$, respectively.

We approximated the zeros until it satisfy the following condition:

\begin{equation} \label{iteration-stop}
\max_{1 \leq k \leq n} \left|P(z_k^{(m)})\right|<10^{-10}.
\end{equation}

In Table I, we give a numerical comparison between several methods
to find all zeros of those polynomials (\ref{test-polynomials}). It contains the
iteration number $m$ and the value $\underset{1 \leq k \leq n}{\max}
\left|P(z_k^{(m)})\right|$ of iterative methods, after we
approximated (\ref{iteration-stop}) to a satisfying label. The smaller the $m$, the
faster approximated on the zeros, When $m$ is the same, it can be
interpreted that a smaller $\underset{1 \leq k \leq n}{\max}
\left|P(z_k^{(m)})\right|$ leads to a higher accuracy of
approximation. All computations have been done using MATLAB.

\begin{center}
Table I. The number of iterations(the error) of iterative methods
\begin{tabular}{l|llllll}
\hline
Poly.  & Method 1 & Method 2 & Method 3 & DFM & NWM  & WLM \\
(\ref{test-polynomials}) & \,\,\, (\ref{method1}) & \,\,\, (\ref{method2}) & \,\,\,(\ref{method3}) & \,\,\,(\ref{derivative-free}) &\,\,\, (\ref{pm}) & \,\,\,(\ref{WDK}) \\
\hline
$P_1$ & 9(8e-14) & 8(9e-14) & 7(1e-14) & 9(9e-14) & 8(1e-10) & 13(3e-12) \\
$P_2$ & 12(7e-13) & 11(4e-13) & 9(4e-13) & 11(2e-13) & 11(1e-12) & 17(2e-12) \\
$P_3$ & 14(2e-11) & 13(5e-12) & 11(8e-12) & 13(8e-11) & 13(7e-12) & 21(2e-11) \\
$P_4$ & 14(2e-11) & 13(2e-11) & 10(2e-11) & 14(2e-11) & 13(2e-11) & 21(2e-11) \\
\hline
\end{tabular}
\end{center}

\section{Conclusion}
In this paper, three new methods for the simultaneous approximation
of all simple zeros of polynomials by utilizing the trapezoidal
Newton's method and the midpoint Newton's method were proposed. It
was proven that each method was of third order. By simultaneously
approximating all simple zeros of polynomials and by comparing
numerical experiments with various methods that are of third order,
we obtained that the results of Method 1 and Method 2 are similar
with that of previous methods. But we found out that the result of
Method 3 are superior than that of any other methods. All methods we
constructed in this paper are new and creative. It seems that these
methods can be applied to various fields, and the study on the
applications of Method 3 is now in progress.

\end{document}